\documentclass[12pt]{article}

\usepackage{amsmath}
\usepackage{mathrsfs}
\usepackage{bbm}
\usepackage{amssymb,a4}
\usepackage{amsfonts,amssymb,amsthm}
\usepackage{extarrows}
\usepackage{color}
\usepackage[all,pdf]{xy}
\usepackage{multicol}
\usepackage{fancyhdr}

\allowdisplaybreaks  

\newtheorem{thm}{\bf Theorem}[section]

\newtheorem{lem}[thm]{\bf Lemma}
\newtheorem{prop}[thm]{\bf Proposition}

\title{\bf\Large Cuspidal modules for solenoidal Lie algebras\\
 over rational quantum tori}
\author{Chengkang Xu$^{\ast}$\\
{\small School of Mathematical Science,
Shangrao Normal University, Shangrao, China
}}

\newcommand{\Z}{{\mathbb Z}}
\newcommand{\C}{{\mathbb C}}
\newcommand{\N}{{\mathbb N}}
\newcommand{\zp}{{\Z_+}}
\newcommand{\g}{\mathfrak{g}}
\newcommand{\U}{{\mathcal U}}
\newcommand{\al}{\alpha}
\newcommand{\be}{\beta}
\newcommand{\dt}{\delta}

\newcommand{\gm}{\gamma}

\newcommand{\e}{{\epsilon}}
\newcommand{\p}{\partial}
\newcommand{\sgm}{\sigma}

\newcommand{\ov}{\overline}
\newcommand{\spanc}[1]{\mathrm{span}_\C\left\{#1\right\}}

\newcommand{\nd}{\N^d}
\newcommand{\ndz}{\nd\backslash\{\bff 0\}}

\newcommand{\f}[2]{F_{#1,#2}}
\newcommand{\pz}{p\Z}
\newcommand{\zd}{\Z^d}
\newcommand{\cd}{\C^d}
\newcommand{\cq}{\C_Q}
\newcommand{\G}{\mathcal{G}}
\newcommand{\call}{\mathcal{L}}
\newcommand{\callx}{\mathcal{L}^{(x)}}
\newcommand{\callt}{\mathcal{L}^{(t)}}
\newcommand{\vabw}{\mathcal V(\al,\be,W)}
\newcommand{\ZZ}{\mathcal{Z}}
\newcommand{\zg}{\ZZ\g}
\newcommand{\gr}{\g_R}
\newcommand{\dercq}{\mathbf{\mathrm{Der}(\cq)}}
\newcommand{\mk}[1]{M_{k_{#1}}(\C)}
\newcommand{\gln}{\mathfrak{gl}_N}
\newcommand{\gld}{\mathfrak{gl}_d}
\newcommand{\gldgm}{\mathfrak{gl}_d(\gm)}
\newcommand{\bff}[1]{{\bf #1}}
\newcommand{\lone}[1]{L_{\bf #1}}             
\newcommand{\ltwo}[2]{L_{{\bf #1}+{\bf #2}}}  
\newcommand{\sgmf}[2]{\sgm({\bf #1},{\bf #2})}
\newcommand{\parder}[1]{\frac{\p}{\p #1}}     
\newcommand{\tone}[1]{t^{\bf #1}}    
\newcommand{\ttwo}[2]{t^{{\bf #1}+{\bf #2}}}
\newcommand{\tthree}[3]{t^{{\bf #1}+{\bf #2}+{\bf #3}}}
\newcommand{\tfour}[4]{t^{{\bf #1}+{\bf #2}+{\bf #3}+{\bf #4}}}
\newcommand{\ovbff}[1]{\ov{\bff #1}}
\newcommand{\tonebar}[1]{t^{\ovbff #1}}
\newcommand{\ttwobar}[2]{t^{{\ovbff #1}+{\bf\ov #2}}}

\newcommand{\xone}[1]{x^{\bf #1}}    
\newcommand{\xtwo}[2]{x^{{\bf #1}+{\bf #2}}}

\newcommand{\bxone}[1]{X^{\bf #1}}    
\newcommand{\bxtwo}[2]{X^{{\bf #1}+{\bf #2}}}

\newcommand{\xonedgm}[1]{\xone #1 d_\gm}

\newcommand{\done}[1]{D({\bff #1})}
\newcommand{\dtwo}[2]{D({\bff #1},{\bff #2})}

\newcommand{\ptwo}[2]{{\p_{\bff #1}^{\bff #2}}}

\newcommand{\inr}{\in R}
\newcommand{\notinr}{\notin R}

\newcommand{\pf}[1]{\noindent{\bf Proof.}#1\hfill{}$\Box$}

\date{}

\begin{document}
\maketitle
\renewcommand{\thefootnote}

\setcounter{footnote}{-1}\footnote{* Corresponding author;
  \emph{E-mail:} xiaoxiongxu@126.com}

{\noindent\bf \normalsize Abstract:}  {\small
In this paper we classify all irreducible cuspidal modules
over a solenoidal Lie algebra over a rational quantum torus,
generalizing the results in \cite{BF2, Su} and \cite{Xu2}.
}

\medskip
{\noindent \small {\bf Keywords}: }{\small solenoidal Lie algebra; cuspidal module;
quantum torus; module of tensor fields.}


\section{Introduction}
\def\theequation{1.\arabic{equation}}
\setcounter{equation}{0}

The classification of irreducible Harish-Chandra modules for infinite dimensional Lie algebras
is a classical problem in the representation theory,
and was solved for the Virasoro algebra $Vir$\cite{M}
and its various generalizations,
such as the twisted Heisenberg-Virasoro algebra\cite{LZ2},
the gap-$p$ Virasoro algebras $Vir_p$\cite{Xu2},
the Witt algebras $W_d$\cite{BF1},
the higher rank Virasoro algebras(or solenoidal Lie algebras) $W_\mu$\cite{LZ1},
and so on.

In order to achieve the classification result for $W_d$ in \cite{BF1},
the authors applied a new technique, called cover method.
Using this cover method, they also classified all irreducible cuspidal modules for $W_\mu$ in \cite{BF2},
which was originally done in \cite{Su}, but in a much more computational way.
Let $A=\C[x_1^{\pm1},\cdots,x_d^{\pm 1}]$ be the Laurent polynomial ring with $d$ variables,
where $d$ is a positive integer,
and let $\mathcal{G}$ stand for $W_d$ or $W_\mu$.
An $A\G$-module is defined to be a module over the Lie algebra $\G\ltimes A$
with an associative $A$-action.
Then the $A$-cover $\hat M$ of a $\G$-module $M$ is a particular $A\G$-quotient
of the module $\G\otimes M$.
Here $A$ plays the role of a coordinate algebra.
It turns out that $\hat M$ is cuspidal if $M$ is,
and $M$ becomes a $\G$-quotient of $\hat M$ provided that $M$ is irreducible.
This reduces the classification of irreducible cuspidal $\G$-modules to
the classification of irreducible cuspidal $A\G$-modules.
It is well known that irreducible cuspidal $AW_d$-modules
are the so-called modules of tensor fields,
originally constructed by Shen\cite{S}, Larsson\cite{Lar},
and further studied in \cite{R}.
While for $W_\mu$, it was proved in \cite{BF2} that cuspidal $AW_\mu$-modules
with the support being one coset of $\zd$ in $\cd$
are in a one-to-one correspondence to finite dimensional modules
over an subalgebra of $\text{Der}\C[x_1,\cdots,x_d]$.
These finite dimensional modules must be one dimensional
if the corresponding $AW_\mu$-module is irreducible,
which makes the irreducible cuspidal $AW_\mu$-module
has weight multiplicities no more than one.

In the present paper we use a modified cover method to classify
all irreducible cuspidal modules for a solenoidal Lie algebra $\g$
over a rational quantum torus $\cq=\C[t_1^{\pm1},\cdots, t_d^{\pm1}]$,
which is a subalgebra of the algebra $\dercq$ of derivations over $\cq$.
The universal central extension of $\g$ was computed in \cite{Xu1}.
When $d=1$, the algebra $\g$ is exactly the gap-$p$ Virasoro algebra studied in \cite{Xu2}.
When constructing the cover for a cuspidal module $M$ for $\g$,
we choose the centre $\ZZ$ of $\cq$, instead of the whole $\cq$,
as the coordinate algebra
and define the $\zg$-cover of $M$ to be
some $\zg$-quotient of the module $\gr'\otimes M$,
where $\gr'$ is an ideal of $\g$.
This construction is a bit different with the $W_\mu$ or the $W_d$ case.

We establish a one-to-one correspondence between
cuspidal modules over the solenoidal Lie algebra $\g$
with support lying in one coset of $\zd$ in $\cd$
and finite dimensional $\Gamma$-graded modules
over an subquotient algebra $\call$ of
the Lie algebra $\text{Der}\C[x_1,\cdots,x_d,t_1,\cdots,t_d]$,
where $\Gamma$ is a finite group closely related to $\cq$.
Through this correspondence we show that
any irreducible cuspidal $\zg$-module
is isomorphic to the module of tensor fields $\vabw$,
for some $\al\in\cd, \be\in\C$ and
finite dimensional $\Gamma$-graded irreducible $\gln$-module $W$,
where $N$ is the order of $\Gamma$.
These modules $\vabw$ may have weight multiplicities larger than one.
The Lie algebra $\call$ has a quotient isomorphic to the direct sum of $\gln$
and a solvable subalgebra of $\gld$.
This is how the seemingly peculiar algebra $\gln$ comes into the picture.
Furthermore, using the modified cover method we prove that
any irreducible cuspidal $\g$-module is isomorphic to
the unique irreducible subquotient of some $\vabw$.

Since the relation between the algebra $\g$ and $\dercq$ is similar to
that between the solenoidal Lie algebra $W_\mu$ and the Witt algebra $W_d$,
we believe our result may play a role in classification of
irreducible Harish-Chandra $\dercq$-modules.
We mention that in \cite{LiuZ} irreducible cuspidal $\dercq$-modules were classified,
also using a cover method.

The paper is arranged as follows.
In Section 2 we recall some results about the gap-$p$ Virasoro algebra $Vir_p$,
the solenoidal Lie algebra $W_\mu$, $\cq$ and $\gln$.
In section 3 we give the construction of modules of tensor fields $\vabw$ over $\g$,
and their irreducibility criterion.
Section 4 is devoted to finite dimensional $\Gamma$-graded $\call$-modules.
In Section 5 we classify all irreducible cuspidal $\zg$-modules,
and in the last section all irreducible cuspidal $\g$-modules are classified
using the cover method.

Throughout this paper, $\C,\Z,\N,\zp$ refer to the set of complex numbers,
integers, nonnegative integers and positive integers respectively.
For a Lie algebra $\mathcal{G}$, we denote by $\U(\mathcal G)$
the universal enveloping algebra of $\mathcal{G}$.
Let $d\in \zp$ and we fix a standard basis $\e_1,\cdots,\e_d$ for the space $\cd$.
Denote by $(\cdot\mid \cdot)$ the inner product on $\cd$.

\section{Notations and Preliminaries}
\def\theequation{2.\arabic{equation}}
\setcounter{equation}{0}

For a Lie algebra,
a weight module is called a Harish-Chandra module
if all weight spaces are finite dimensional,
and called a {\bf cuspidal} module if all weight spaces are uniformly bounded.
The {\bf support} of a weight module is defined to be the set of all weights.

\subsection{The gap-$p$ Virasoro algebra $Vir_p$}

Let $p$ be a positive integer
and $\C[x^{\pm1}]$ denote the Laurent polynomial ring in one variable $x$.
The gap-$p$ Virasoro algebra $Vir_p$ is a Lie algebra with a basis
$$\{x^{m+1}\parder x,x^s, C_i\mid m\in\pz,s\notin\pz, 0\leq i\leq p-1\},$$
and Lie brackets
$$\begin{aligned}
  &[x^{m+1}\parder x,x^{n+1}\parder x]
     =(n-m)x^{m+n+1}\parder x+\dt_{m+n,0}\frac{1}{12}
       \left((\frac{m}{p})^3-(\frac{m}{p})\right)C_0;\\
  &[x^{m+1}\parder x,x^r]=rx^{m+r};\ \ \ [x^r,x^s]=\dt_{r+s,0}rC_{\overline r},
\end{aligned}$$
where $m,n\in \pz, r,s\notin \pz$ and
we use ${\overline r}$ to represent the residue of $r$ by $p$.

We recall from \cite{Xu2} the module of intermediate series over $Vir_p$.
Let $F=(\f ij)$ be a $(p-1)\times p$ complex matrix,
with index $1\leq i\leq p-1,\ 0\leq j\leq p-1$,
satisfying the following three conditions
\begin{description}
 \item[(I)] $0\in o(F)=\{j\mid \f ij\neq 0\text{ for some } i\}$;
 \item[(II)] if $\f ij\neq 0$ then $\f s {\ov{i+j}}\neq 0$ for some $1\leq s\leq p-1$;
 \item[(III)] $\f r{\ov{i+s}}\f si=\f s{\ov{i+r}}\f ri$
                for any $0\leq i\leq p-1$ and $1\leq r,s\leq p-1$;
\end{description}
For any $j\in o(F)$ denote the space $V_{(j)}=\spanc{v_{j+pk}\mid k\in\Z}$.
Let $a,b\in \C$ and
define the $Vir_p$-module structure on $\bigoplus_{j\in o(F)}V_{(j)}$ by
$$\begin{aligned}
  &x^{m+1}\parder x v_{j+n}=(a+j+n+mb)v_{j+n+m};\\
  &x^sv_{j+n}=\f{\ov s}jv_{j+n+s};\ \ C_iv_{j+n}=0,
\end{aligned}$$
where $m,n\in \pz,\ s\notin \pz, j\in o(F)$ and $i=0,1,\cdots,p-1$.
We denote this module by $V(a,b,F)$, and call it a module of intermediate series over $Vir_p$.
Roughly speaking, the module $V(a,b,F)$ is a sum of several modules of intermediate series
over the subalgebra $\spanc{x^{m+1}\parder x,C_0\mid m\in\pz}$,
which is isomorphic to the normal Virasoro algebra.
Let $P=\spanc{x^m\mid m\in\pz}$.
In \cite{Xu2}, it was proved that
any irreducible cuspidal modules over $Vir_p\ltimes P$
with an associative $P$-action must be of the form $V(a,b,F)$.

\subsection{The solenoidal Lie algebra $W_\mu$}

Let $d\in \zp$ and $\mu=(\mu_1,\cdots,\mu_d)^T\in\cd$ is called {\bf generic}
if $\mu_1,\cdots,\mu_d$ are linearly independent over the field of rational numbers.
Let $A=\C[x_1^{\pm1},\cdots,x_d^{\pm 1}]$ be the Laurent polynomial ring with $d$ variables.
For any $\bff m=(m_1,\cdots,m_d)^T\in\zd$ denote $\xone m=x_1^{m_1}\cdots x_d^{m_d}$.
Set $\p_x=\sum_{i=1}^d\mu_ix_i\parder{x_i}$.
Then the solenoidal Lie algebra $W_\mu$ over $A$ has a basis
$$\{\xone m\p_x\mid\bff m\in\zd\}$$
and Lie bracket
$$[\xone m\p_x,\xone n\p_x]=(\mu\mid\bff n-\bff m)\xtwo mn\p_x.$$
Let $\al\in\cd,\be\in\C$.
The modules $T(\al,\be)$ of tensor fields over $W_\mu$ have bases
$\{v_{\bff s}\mid \bff s\in\zd\}$ and the $W_\mu$-action
$$(\xone m\p_x)\cdot v_{\bff s}=(\mu\mid \al+\bff s+\be\bff m)v_{\bff m+\bff s}.$$
It is well known that the module $T(\al,\be)$ is reducible if and only if
$\al\in\zd$ and $\be\in\{0,1\}$.
For $\al\in\zd$, $T(\al,0)$ has a unique irreducible quotient $T(\al,0)/\C v_{-\al}$,
and $T(\al,1)$ has a unique submodule $\spanc{v_{\bff s}\mid \bff s\neq -\al}$
of codimension one.
Moreover, the module $T(\al,\be)$ may be equipped with an additional associative $A$-action
$\xone m v_{\bff s}=v_{\bff m+\bff s}$.

\begin{thm}[\cite{Su,BF2}]
Any irreducible cuspidal module over $W_\mu\ltimes A$ with an associative $A$-action
must be of the form $T(\al,\be)$ for some $\al\in\cd$ and $\be\in\C$.
\end{thm}

\subsection{The quantum torus $\cq$ and the Lie algebra $\gln$}

Let $Q=(q_{ij})$ be a $d\times d$ complex matrix
with all $q_{ij}$ being roots of unity and satisfying
$$q_{ii}=1,\ \ \ q_{ij}q_{ji}=1\text{ for all }1\leq i,j\leq d.$$
The rational quantum torus relative to $Q$ is the unital associative algebra
$\cq=\C[t_1^{\pm1},\cdots, t_d^{\pm1}]$
with multiplication
$$t_it_j=q_{ij}t_jt_i\text{ for all }1\leq i,j\leq d.$$
For an element ${\bff m}=(m_1,\cdots, m_d)^T\in\zd$ we denote $\tone m=t_1^{m_1}\cdots t_d^{m_d}$.

For $\bff m, \bff n \in \zd$, denote
$$\sgmf mn=\prod_{1\leq i<j\leq d}q_{ji}^{n_jm_i}\text{ and }
R=\{\bff m\zd\mid \sgmf mn=\sgmf nm\text{ for any }\bff n\in\zd\}.$$
Clearly, the center $\ZZ$ of $\cq$ is spanned by $\{\tone m\mid\bff m\in R\}$.
From Theorem 4.5 in \cite{N}, up to an isomorphism of $\cq$,
we may always assume that
$q_{2i,2i-1}=q_i, q_{2i-,2i}=q_i^{-1}$ for $1\leq i\leq z$,
and other entries of $Q$ are all 1,
where $z\in\zp$ with $2z\leq d$
and the orders $k_i$ of $q_i$ as roots of unity satisfy
$k_{i+1}\mid k_i$ for $1\leq i\leq z$.
Then the subgroup $R$ of $\zd$ has a simple form
$$R=\bigoplus_{i=1}^z\left(\Z k_i\e_{2i-1}\oplus \Z k_i\e_{2i}\right)
\oplus \bigoplus_{l>2z}\Z\e_l.$$
Moreover, we have
$$\sgmf mr=\sgmf rm,\ \ \tone m\tone r=\ttwo mr\text{ for all }\bff m\in R, \bff r\in\zd.$$

Let $\mathcal I=\spanc{\ttwo nr-\tone r\mid n\in R,r\in\zd}$,
which is an ideal of the associative algebra $\cq$.
By \cite{N} and \cite{Z}, we have
$\cq/\mathcal I\cong \bigotimes_{i=1}^z\mk i\cong M_N(\C)$,
where $N=\prod_{i=1}^zk_i$ and
$M_n(\C)$ is the associative algebra of all $n\times n$ complex matrices.
It is well known that $\mk i$ can be generated by
$$\begin{aligned}
&X_{2i-1}=E_{1,1}+q_iE_{2,2}+\cdots+q_i^{k_1-1}E_{k_i,k_i},\\
&X_{2i}=E_{1,2}+E_{2,3}+\cdots+ E_{k_i-1,k_i}+E_{k_i,1},
\end{aligned}$$
where $E_{kl}$ represents the $k_i\times k_i$ matrix with 1 in the $(k,l)$-entry
and 0 elsewhere.
Let $E$ denote the identity matrix of suitable order.
It is easy to see that $X_{2i}^{k_i}=X_{2i-1}^{k_i}=E$
and $X_{2i}X_{2i-1}=q_i X_{2i-1}X_{2i}$.
Denote $X^{\bff n}=\bigotimes_{i=1}^z X_{2i-1}^{n_{2i-1}}X_{2i}^{n_{2i}}$
for $\bff n\in\zd$.
Notice that for any $\bff n\in R$, $\bxone n$ is the identity matrix in $M_N(\C)$,
and $\bxone r\bxone s=\sgmf rs\bxtwo rs$.

Define the Lie bracket $[a,b]=ab-ba$ on $M_N(\C)$,
and we write $\gln$ for $M_N(\C)$ as a Lie algebra.
We have the Lie bracket for $\gln$
$$[\bxone m,\bxone n]=(\sgmf mn-\sgmf nm)\bxtwo mn.$$
Set $\Gamma=\zd/R$ and let $\ovbff n$ denote the image of $\bff n$.
Clearly, $\Gamma$ has order $N$,
and $\gln$ is a $\Gamma$-graded Lie algebra with homogeneous spaces
$(\gln)_{\ovbff n}=\spanc{\bxone n}$.
In this paper by a $\Gamma$-gradation on $\gln$ we always mean this gradation.
Set
$$\Gamma_0=\{\bff n\in\zd\mid 0\leq n_{2i-1}<k_i,0\leq n_{2i}<k_i,1\leq i\leq z,
\text{ and }n_l=0,2z<l\leq d\},$$
which is a complete set of representatives for $\Gamma$.
When a representative of $\ovbff n$ is needed we always choose
the one $\bff n\in\Gamma_0$.

\section{Modules of tensor fields for $\g$}
\def\theequation{3.\arabic{equation}}
\setcounter{equation}{0}

In this section we construct modules of tensor fields $\vabw$ for
the solenoidal Lie algebra $\g$ over the quantum torus $\cq$,
and give an irreducibility criterion for $\vabw$.

Let $\gm=(\gm_1,\cdots,\gm_d)^T\in\cd$ be generic and set
$$\lone m=\begin{cases}
\tone m\sum_{i=1}^d\gm_it_i\parder{t_i} &  \text{ if }m\in R;\\
     \tone m                            &  \text{ if }m\notin R,
\end{cases}$$
where $\tone m$ is the inner derivation on $\cq$ defined by
$\tone m(\tone r)=\sgmf mr\ttwo mr$.
The Lie algebra $\g$ we consider in this paper has a basis
$\{\lone m\mid \bff m\in\zd\}$,
subject to the Lie brackets
$$\begin{aligned}
&[\lone m,\lone n]=(\gm\mid{\bff n}-{\bff m}) \ltwo mn;\\
&[\lone m,\lone s]=(\gm\mid\bff s) \ltwo ms;\ \
[\lone r,\lone s]=(\sgmf rs-\sgmf sr) \ltwo rs,
\end{aligned}$$
for $\bff m,\bff n\in R, \bff r, \bff s\notin R$.
The subalgebra $\gr$ of $\g$ spanned by $\{\lone m\mid \bff m\inr\}$
is isomorphic to the solenoidal Lie algebra $W_{\mu}$,
where $\mu=B\gm, B=\mathrm{diag}\{k_1,k_1,k_2,k_2,\cdots,k_z,k_z,1,\cdots,1\}$.
This isomorphism is given by $\lone m\mapsto \xone n\sum_{i=1}^d \mu_ix_i\parder{x_i}$,
where $\bff m=B\bff n$.
The Lie algebra $\g$ may be considered as a quantum version of $W_\mu$.

Now we give the construction of the module of tensor fields for $\g$.
Let $\al\in\cd, \be\in\C$ and $W=\bigoplus_{\bff s\in\Gamma}W_{\bff s}$
be a $\Gamma$-graded $\gln$-module.
Recall the center $\ZZ$ of $\cq$.
Define a $\g$-module structure on
$$\vabw=\bigoplus_{\ovbff s\in\Gamma}W_{\ovbff s}\otimes \tone s\ZZ$$
by(here $\bff s$ is the representative of $\ovbff s$ in $\Gamma_0$)
$$ \begin{aligned}
   &\lone m(w_{\ovbff s}\otimes\ttwo ns)
      =(\gm\mid\al+\bff n+\bff s+\be\bff m)w_{\ovbff s}\otimes\tthree mns;\\
   &\lone r(w_{\ovbff s}\otimes\ttwo ns)=(\bxone rw_{\ovbff s})\otimes\tthree rns,
 \end{aligned}$$
where $\bff m,\bff n\inr, \bff r\notinr, \bff s\in\zd$ and $w_{\ovbff s}\in W_{\ovbff s}$.
We call $\vabw$ a module of tensor fields over $\g$.

\begin{lem}\label{lem3.1}
 If the $\Gamma$-graded $\gln$-module is irreducible and $\dim W>1$,
 then the $\g$-module $\vabw$ is irreducible for any $\al\in\cd, \be\in\C$.
\end{lem}
\pf{
Let $U$ be a nonzero $\g$-submodule of $\vabw$ and
$0\neq w_{\ovbff s}\otimes \tone s\in U$ for some $w_{\ovbff s}\in W_{\ovbff s}, \bff s\in\zd$.
We claim that
$$w_{\ovbff k}\otimes\ttwo km\in U\text{ for any }
  \bff s\neq \bff k\in\zd\text{ and }\bff m\inr.$$
Since $W$ is a $\Gamma$-graded irreducible $\gln$-module,
we see that $w_{\ovbff k}=aw_{\ovbff s}$
for some homogeneous element $a$ in $\U(\gln)$ of the form
$$a=\sum c_{{\bff s_1}\cdots {\bff s_p}}X^{\bff s_1}\cdots X^{\bff s_p},$$
where $c_{{\bff s_1}\cdots {\bff s_p}}\in\C$ and
the sum takes over finitely many $(\bff s_1,\cdots,\bff s_p)$
such that $\bff s_1+\cdots\bff s_p=\bff k-\bff s$.
Denote $b=\sum c_{{\bff s_1}\cdots {\bff s_p}}L_{\bff s_1}\cdots L_{\bff s_1}\in\U(\g)$.
Then we have
$$w_{\ovbff k}\otimes\tone k=aw_{\ovbff s}\otimes\tone k=b(w_{\ovbff s}\otimes\tone s).$$
Since $\lone m(w_{\ovbff k}\otimes\tone k)=
 (\gm\mid \al+\bff k+\be\bff m)(w_{\ovbff k}\otimes\ttwo mk)\in U$,
we get that
$$w_{\ovbff k}\otimes\ttwo mk\in U\text{ if }\al+\bff k+\be\bff m\neq\bff 0.$$
Assume $\al+\bff k+\be\bff m=\bff 0$.
Since
$$\lone m(w_{\ovbff s}\otimes\tone s)=
 (\gm\mid \al+\bff s+\be\bff m)(w_{\ovbff s}\otimes\ttwo ms)=
 (\gm\mid \bff s-\bff k)(w_{\ovbff s}\otimes\ttwo ms)\in U$$
we obtain $w_{\ovbff s}\otimes\ttwo ms\in U$.
Then we still have $w_{\ovbff k}\otimes\ttwo mk=b(w_{\ovbff s}\otimes\ttwo ms)\in U$.

We still need to show $w'_{\ovbff s}\otimes\ttwo ms\in U$ for any $\bff m\inr$
and $w'_{\ovbff s}\in W_{\ovbff s}$.
But this is just a replicate of the proof of the above claim starting with
a vector $w_{\ovbff k}\otimes\tone k\in U$ with $\bff k\neq \bff s$.
So $U=\vabw$ and the lemma stands.
}

Furthermore, we have the following
\begin{thm}
  Let $\al\in\cd,\be\in\C$ and $W$ be a $\Gamma$-graded irreducible $\gln$-module.
  The $\g$-module $\vabw$ is reducible
  if and only if $\dim W=1, \al\in\zd$ and $\be\in\{0,1\}$.
\end{thm}
\pf{
By Lemma \ref{lem3.1} we only need consider when $\dim W=1$.
In this case $\lone s\vabw=0$ for any $\bff s\notinr$,
and hence $\vabw$ reduces to a module over the subalgebra $\gr$ of $\g$,
which is isomorphic to the solenoidal Lie algebra $W_{B\gm}$.
Then the theorem follows from the irreducibility criterion of the $W_{B\gm}$-module $T(\al,\be)$.
}

\section{The algebra $\call$ and finite dimensional modules}
\def\theequation{4.\arabic{equation}}
\setcounter{equation}{0}

In this section we introduce the related algebra $\call$
and study finite dimensional modules over $\call$.

Consider the associative algebra $AC=\C[x_1,\cdots,x_d]\oplus(\cq/\ZZ)$
with $x_it_j=t_jx_i$ for all $1\leq i,j\leq d$.
We simply write $\tonebar s$ for the image of $\tone s$ in $\cq/\ZZ$.
Denote by $\call$ the subalgebra of $\mathrm{Der}(AC)$
spanned by
$$\{\xonedgm m,\xone n\tonebar s\mid \bff m\in \ndz, \bff n\in \nd, \ovbff s\in\Gamma\},$$
where $d_\gm=\sum_{i=1}^d \gm_i\parder{x_i}+\sum_{i=1}^d \gm_it_i\parder{t_i}$.
The Lie bracket of $\call$ is
\begin{equation}\label{eq4.2}
 \begin{aligned}
   &[\xonedgm m,\xonedgm n]=\sum_{i=1}^d\gm_i(n_i-m_i)x^{\bff m+\bff n-\e_i}d_\gm;\\
   &[\xonedgm m,\xone l\tonebar s]=\sum_{i=1}^d\gm_il_ix^{\bff m+\bff l-\e_i}\tonebar s+
                                 (\gm\mid\bff s)\xtwo ml\tonebar s;\\
   &[\xone p\tonebar r,\xone l\tonebar s]=(\sgmf rs-\sgmf sr)\xtwo pl\ttwobar rs,
 \end{aligned}
\end{equation}
where $\bff m,\bff n\in\ndz,\bff p,\bff l\in\nd$ and $\ovbff r,\ovbff s\in\Gamma$.
Here we have chosen the representative for $\ovbff s\in\Gamma$ in $\Gamma_0$ as usual.
The algebra $\call$ has a $\Gamma$-gradation
$\call=\bigoplus_{\ovbff s\in\Gamma}\call(\ovbff s)$ where
$$\begin{aligned}
   &\call(\ovbff s)=\spanc{\xone n\tonebar s\mid \bff n\in\nd};\\
   &\call(\ovbff 0)=\spanc{\xonedgm m,\xone n\tonebar 0\mid \bff m\in\ndz,\bff n\in\nd}.
 \end{aligned}$$
Denote two subalgebras of $\call$
$$
\callx=\spanc{\xonedgm m\mid \bff m\in\ndz},\ \
\callt=\spanc{\xone n\tonebar s\mid \bff n\in\nd,\ovbff s\in\Gamma}.
$$
Set $\deg(x_i)=1,\deg(d_\gm)=-1$ and $|\bff m|=\sum_{j=1}^dm_j$
for all $1\leq i\leq d$ and $\bff m\in\nd$.
We get $\Z$-gradations for the algebras $\callx=\bigoplus_{i\geq 0}\callx_i$
and $\callt=\bigoplus_{i\geq 0}\callt_i$ such that
$$
\callx_i=\spanc{\xonedgm m\mid |\bff m|=i+1},\ \
\callt_i=\spanc{\xone n\tonebar s\mid |\bff n|=i,\ovbff s\in\Gamma}.
$$
Set $\call_i=\callx_i\oplus\callt_i$.
Notice that the derivation $d_\gm$ is not homogeneous,
hence the algebra $\call$ is not $\Z$-graded.
But the subspace $\call_+=\bigoplus_{i\geq 1}\call_i$ still makes an ideal of $\call$.
The main result in this section is the following
\begin{thm}\label{thm4.1}
(1) The commutator
    $[\callx,\callx]=[\callx_0,\callx_0]\oplus \left(\bigoplus_{j\geq 1}\callx_j\right)$;\\
(2) Every finite dimensional representation $(U,\rho)$ for $\callx$
    satisfies $\rho(\callx_p)=0$ for $p\gg0$;\\
(3) Every finite dimensional irreducible $\callx$-module is one dimensional,
    and parametrized by some $\be\in\C$ such that $x_i\parder{x_i}\mapsto \be\gm_i, 1\leq i\leq d$
    and $\callx_j\mapsto 0$ for $j\geq 1$;\\
(4) Every finite dimensional representation $(U,\rho)$ for $\call$
    satisfies $\rho(\call_p)=0$ for $p\gg0$;\\
(5) The ideal $\call_+$ annihilates every finite dimensional irreducible $\call$-module.
\end{thm}
\pf{
The first three statements are exactly the Theorem 3.1 from \cite{BF2}.
For (4), consider $U$ as a $\callx$-module.
By (3) we have $\rho(\callx_p)=0$ for $p\gg0$.
Then it follows from $[\xonedgm m,\tonebar s]=(\gm\mid\bff s)\xone m\tonebar s$
that $\rho(\xone m\tonebar s)=0$ if $|\bff m|>p$.
This proves (4).

Let $(V,\rho)$ be a finite dimensional irreducible representation of $\call$.
By (4) and equation (\ref{eq4.2}) we see that
$\rho(\call_+)$ is a finite dimensional nilpotent Lie algebra.
Let $V'=\{v\in V\mid \call_+v=0\}$,
which is a $\call$-submodule of $V$.
Moreover, by Lie's Theorem, there exists a common eigenvector $v\in V$ for $\rho(\call_+)$
such that $\rho(a)v=0$ for all $a\in\call_+$.
This means $V'\neq0$.
So $V'=V$ by the irreducibility of $V$.
Hence $\call_+V=0$.
}

For later use we mention that there is an isomorphism
from the quotient algebra $\call/\call_+$ to $\gldgm\oplus\gln$,
where $\gldgm=\spanc{\e_i\gm^T\in\gld\mid 1\leq i\leq d}$,
defined by
\begin{equation}\label{eq4.2}
  x_id_\gm+\call_+\mapsto \e_i\gm^T,\ \ \ \tonebar s+\call_+\mapsto \bxone s.
\end{equation}
Notice that $\gldgm$ is solvable since $[\gldgm,\gldgm]$ is abelian.

\section{Cuspidal $\zg$-modules}
\def\theequation{5.\arabic{equation}}
\setcounter{equation}{0}

In this section we study a specific class of cuspidal modules over $\g$,
and prove that the irreducible ones are exactly those modules of tensor fields.

Recall the centre $\ZZ$ of $\cq$.
We may form an extended Lie algebra $\g\ltimes\ZZ$.
We call a module $V$ for $\g\ltimes\ZZ$ a $\zg$-module
provided that the $\ZZ$-action on $V$ is associative.

The following are two examples of $\zg$-module.
Set $\gr'=\spanc{\lone s\mid \bff s\notinr}$,
which is an ideal of $\g$.
Then $\gr'$ makes a $\zg$-module if we define
$$\lone m\cdot\lone s=[\lone m,\lone s],\ \ \ \tone n\cdot\lone s=\ltwo ns\ \ \
\text{for }\bff m\in\zd,\bff n\inr\text{ and }\bff s\notinr.$$
The second example is the module
$\vabw=\bigoplus_{\ovbff s\in\Gamma}W_{\ovbff s}\otimes\tone s\ZZ$
of tensor fields for $\g$ with a $\ZZ$-action given by
$$\tone m(w_{\ovbff s}\otimes\ttwo ns)=w_{\ovbff s}\otimes\tthree mns,$$
for $\bff m,\bff n\inr,\bff s\notinr$ and $w_{\ovbff s}\in W_{\ovbff s}$.
Clearly, if $W$ is a $\Gamma$-graded irreducible $\gln$-module,
then $\vabw$ is irreducible as a $\zg$-module for any $\al\in\cd,\be\in\C$.

Now let $M$ be a cuspidal $\zg$-module with support lying in $\al+\zd$.
Since the $\ZZ$-action is associative, $M$ can be represented as
$$M\cong \bigoplus_{\ovbff s\in\Gamma}U_{\ovbff s}\otimes\ZZ,$$
where $U_{\ovbff s}=M_{\al+\bff s}$ for $\bff s\in\Gamma_0$.
Set $U=\bigoplus_{\ovbff s\in\Gamma}U_{\ovbff s}$.

For later use we should need some operators on $M$.
For $\bff m\inr,\bff r\notinr$, consider
$$\done m=t^{-\bff m}\lone m,\ \ \dtwo mr=t^{-\bff m}\ltwo mr.$$
The operator $\done m$ may be restricted to each $U_{\ovbff s}$,
and $\dtwo mr$ may be restricted to an operator
$U_{\ovbff s}\longrightarrow U_{\ovbff {r+s}}$.
The operators $\done m,\dtwo mr$ completely determine the $\g$-action on $M$, since
$$\begin{aligned}
  &\lone m(\tone n v_{\ovbff s})
    =(\gm\mid \bff n)\ttwo mn v_{\ovbff s}+\ttwo mn\done mv_{\ovbff s}\\
  &\ltwo mr(\tone n v_{\ovbff s})
    =\tone n(\ltwo mr v_{\ovbff s})=\ttwo mn(\dtwo mrv_{\ovbff s}).
  \end{aligned}$$
The following lemma is easy to check.

\begin{lem}\label{lem5.1}
 For $\bff m,\bff n\inr$ and $\bff r,\bff s\notinr$, we have
 $$\begin{aligned}
  &[\done m,\done n]=(\gm\mid\bff m)(D(\bff m+\bff n)
                    -\done m)-(\gm\mid\bff n)(D(\bff m+\bff n)-\done n);\\
  &[\done m,\dtwo ns]=(\gm\mid\bff n+\bff s)D(\bff m+\bff n,\bff s)-(\gm\mid\bff n)\dtwo ns;\\
  &[\dtwo mr,\dtwo ns]=(\sgmf rs-\sgmf sr)D(\bff m+\bff n,\bff r+\bff s).
 \end{aligned}$$
\end{lem}

\begin{prop}\label{prop5.2}
 Let $M$ be a cuspidal $\zg$-module with support lying in $\al+\zd$.
 We may write $M\cong \bigoplus_{\ovbff s\in\Gamma}U_{\ovbff s}\otimes\ZZ,$
 where $U_{\ovbff s}=M_{\al+\bff s}$ for $\bff s\in\Gamma_0$.
 Then the action of $\g$ on $M$ is given by
 $$\begin{aligned}
   &\lone m(v_{\ovbff s}\otimes \tone n)
         =((\gm\mid\bff n)+\done m)v_{\ovbff s}\otimes \ttwo mn;\\
   &\ltwo mr(v_{\ovbff s}\otimes \tone n)=\dtwo mr\cdot v_{\ovbff s}\otimes \ttwo mn,
 \end{aligned}$$
 for $\bff m,\bff n\inr, \bff r\in\Gamma_0\backslash\{\bff 0\}$ and $v_{\ovbff s}\in U_{\ovbff s}$.
 Here the operators $\done m: U_{\ovbff s}\longrightarrow U_{\ovbff s}$
 for each $\ovbff s\in\Gamma$, can be expressed as
 $\mathrm{End}(U_{\ovbff s})$-valued polynomial in $\bff m$ with constant term
 $\done 0=(\gm\mid\al+\bff s)Id$,
 and the operators $\dtwo mr: U_{\ovbff s}\longrightarrow U_{\ovbff {r+s}}$
 for each $\ovbff s\in\Gamma$, can be expressed as
 $\mathrm{Hom}(U_{\ovbff s},U_{\ovbff {r+s}})$-valued polynomial in $\bff m$
 whose constant term $\dtwo 0r$ has the form of an upper triangular matrix
 relative to suitable bases of $U_{\ovbff s}$ and $U_{\ovbff {r+s}}$.
\end{prop}
\pf{
Notice that $\gr$ is isomorphic to the solenoidal Lie algebra $W_{B\gm}$.
The part concerning the action of $\gr$ and $\done m$ is just
what was stated in Theorem 4.5 in \cite{BF1}.
Next we prove the $\gr'$-action by induction on $d$.
If $d=1$, then the algebra $\g$ shrinks either to a gap-$p$ Virasoro algebra(for suitable $p$),
in which case the result follows from \cite{Xu2},
or to a normal Virasoro algebra, in which case $\gr'=0$ and there is nothing to prove.

Now let us establish the induction step as follows.
By induction assumption the operators
$D(\bff m-m_i\e_i, \bff r-r_i\e_i)$ have polynomial dependence
on $m_1,\cdots,m_{i-1},m_{i+1},\cdots,m_d$,
and $D(m_i\e_i, r_i\e_i)$ is a polynomial on $m_i$,
for all $i=1,2,\cdots,d$ and $\bff r\notinr,r_i\e_i\notinr$ and $\bff r-r_i\e_i\notinr$.
Note that $\bff m,m_i\e_i\inr$.
From
$$ [D(\bff 0,r_i\e_i),D(\bff m-m_i\e_i,\bff r-r_i\e_i)]
  =(\sgm(r_i\e_i,\bff r-r_i\e_i)-\sgm(\bff r-r_i\e_i, r_i\e_i))D(\bff m-m_i\e_i, \bff r)
$$
we see that $D(\bff m-m_i\e_i, \bff r)$ is a polynomial
on $m_1,\cdots,m_{i-1},m_{i+1},\cdots,m_d$.
Then consider
$$
 [D(m_i\e_i),D(\bff m-m_i\e_i,\bff r)]=(\gm\mid\bff m-m_i\e_i+\bff r)D(\bff m, \bff r)
  -(\gm\mid\bff m-m_i\e_i)D(\bff m-m_i\e_i, \bff r)
$$
we see that $D(\bff m, \bff r)$ may be expressed as
a rational function $\frac{P_i(\bff m)}{F_i(\bff m)}$, where
$$F_i(\bff m)=(\gm\mid\bff m-m_i\e_i+\bff r)=\gm_ir_i+\sum_{j\neq i}\gm_j(m_j+r_j)\neq 0$$
since $\gm$ is generic and $m_j+r_j\neq 0$ for all $j\neq i$.
Especially we have $\frac{P_1(\bff m)}{F_1(\bff m)}=\frac{P_2(\bff m)}{F_2(\bff m)}$
in $\mathrm{Hom}(U_{\ovbff s},U_{\ovbff {r+s}})\otimes \C(m_1,\cdots,m_d)$.
Hence
$$P_1(\bff m){F_2(\bff m)}={P_2(\bff m)}{F_1(\bff m)}.$$
Notice that $F_1(\bff m),F_2(\bff m)$ are coprime to each other and
the ring $\C[m_1,\cdots,m_d]$ is a unique factorization domain.
It follows that $F_1(\bff m)$ divides $P_1(\bff m)$.
So $\dtwo mr$ is a polynomial.

It remains to show that the value of the polynomial $\dtwo mr$ at $\bff m=\bff 0$
coincides with $\dtwo 0r$ on $U_{\ovbff s}$.
Let $k$ be the order of the image $\ovbff r\in\Gamma$ of $\bff r$,
and we have $k\bff r\inr$.
Denote by $\g(\bff r)$ the Lie subalgebra of $\g$ generated by
$\{L_{l\bff r}\mid l\in\Z\},$
which is a gap-$k$ Virasoro algebra.
Consider $M$ as a $\g(\bff r)$-module with composition series
$$0=M_0\subset M_1\subset\cdots\subset M_l=M.$$
Each quotient $M_i/M_{i-1}$ is an irreducible cuspidal $\g(\bff r)$-module,
and by \cite{Xu2} has weight multiplicities no more than one.
It is clear that $l$ is the maximal dimension of the weight spaces of $M$.
Then, according to this composition series, one can choose a basis of $M$,
with respect to which the operator $\dtwo 0r=\lone r$
has the form of upper triangular matrices of order $l$
on all weight spaces of $M$, hence on $U_{\ovbff s}$.
Moreover, from the construction of module of intermediate series from \cite{Xu2},
we see that $D(k\bff r,\bff r)$ has the same action as $\dtwo 0r$ on $U_{\ovbff s}$.
Notice that the constant term $\done 0$ of
$D(\bff m-k\bff r)$ acts as $(\gm\mid\al+\bff s)\mathrm{Id}$ on $U_{\ovbff s}$,
and as $(\gm\mid\al+\bff r+\bff s)\mathrm{Id}$ on $U_{\ovbff {r+s}}$.
By considering the constant term in both sides of
$$[D(\bff m-k\bff r),D(k\bff r,\bff r)]=(k+1)(\gm\mid\bff r)\dtwo mr
   -k(\gm\mid\bff r)D(k\bff r,\bff r),$$
we see that the value of $\dtwo mr$ at $\bff m=\bff 0$
equals to $\dtwo 0r$ on $U_{\ovbff s}$.
}

Define operators $\ptwo 0p\in\mathrm{End}(U_{\ovbff s})$
and $\ptwo rp\in\mathrm{Hom}(U_{\ovbff s},U_{\ovbff {r+s}})$
for $\bff r\in\Gamma_0\backslash\{\bff 0\}$,
by expansions of the polynomials $\done m$ and $\dtwo mr$ in $\bff m$:
$$\done m=\sum_{\bff p\in \nd}\frac{\bff m^{\bff p}}{\bff p!}\ptwo 0p;\ \ \ \ \
   \dtwo mr=\sum_{\bff p\in \nd}\frac{\bff m^{\bff p}}{\bff p!}\ptwo rp,$$
with only finite number of the operators $\ptwo 0p, \ptwo rp$ being nonzero.
Here $\bff p!=p_1!\cdots p_d!$ and
$\bff m^{\bff p}=m_1^{p_1}\cdots m_d^{p_d}$.

Now we expand the commutator $[\done m,\done n]$ in $\bff m,\bff n$,
$$\sum_{\bff p,\bff l\in\nd}\frac{\bff m^{\bff p}}{\bff p!}\frac{\bff n^{\bff l}}{\bff l!}
  [\ptwo 0p,\ptwo 0l]=(\gm\mid \bff m)\sum_{\bff q\in\nd}
  \frac{(\bff m+\bff n)^{\bff q}-\bff m^{\bff q}}{\bff q!}\ptwo 0q-
  (\gm\mid \bff n)\sum_{\bff q\in\nd}
  \frac{(\bff m+\bff n)^{\bff q}-\bff n^{\bff q}}{\bff q!}\ptwo 0q.$$
Comparing the coefficients at
$\frac{\bff m^{\bff p}}{\bff p!}\frac{\bff n^{\bff l}}{\bff l!}$
in both sides, we get the Lie bracket of $\ptwo 0p$
\begin{equation}\label{eq5.1}
 [\ptwo 0p,\ptwo 0l]=
 \begin{cases}
  \sum_{i=1}^d\gm_i(l_i-p_i)\p_{\bff 0}^{\bff p+\bff l-\e_i}&\text{ if }\bff p,\bff l\neq\bff 0;\\
  0                        &\text{ if }\bff p=\bff 0\text{ or }\bff l=\bff 0.
 \end{cases}
\end{equation}
Similarly by expanding the commutators $[\done m,\dtwo ns]$ and $[\dtwo mr,\dtwo ns]$,
then comparing coefficients at
$\frac{\bff m^{\bff p}}{\bff p!}\frac{\bff n^{\bff l}}{\bff l!}$
in both sides, we get the Lie bracket of $[\ptwo 0p,\ptwo sl]$
\begin{equation}\label{eq5.2}
 [\ptwo 0p,\ptwo sl]=
 \begin{cases}
  (\gm\mid\bff s)\p_{\bff s}^{\bff p+\bff l}+\sum_{i=1}^d\gm_il_i\p_{\bff s}^{\bff p+\bff l-\e_i}
                                  &\text{ if }\bff p\neq\bff 0;\\
  (\gm\mid\bff s)\ptwo sl         &\text{ if }\bff p=\bff 0,
 \end{cases}
\end{equation}
and the Lie bracket of $[\ptwo rp,\ptwo sl]$
\begin{equation}\label{eq5.3}
 [\ptwo rp,\ptwo sl]=(\sgmf rs-\sgmf sr)\p_{\bff r+\bff s}^{\bff p+\bff l}.
\end{equation}

Recall the Lie bracket of the algebra $\call$.
The equations (\ref{eq5.1}), (\ref{eq5.2}) and (\ref{eq5.3}) imply that
the operators $\{\ptwo 0p, \ptwo sl\mid\bff p\in\ndz,
 \bff s\in\Gamma_0\backslash\{\bff 0\},\bff l\in\nd \}$
yield a finite dimensional $\Gamma$-graded representation of
the algebra $\call$ on the space $U$.
Denote by $\rho$ this representation map.
Combining with Proposition \ref{prop5.2} we get the following

\begin{thm}\label{thm5.3}
There exists an equivalence between the category of finite dimensional $\Gamma$-graded
$\call$-modules and the category of cuspidal $\zg$-modules with support
lying in some coset $\al+\zd$.
The equivalence functor associates to a finite dimensional $\Gamma$-graded $\call$-module
$U=\bigoplus_{\ovbff s\in\Gamma}U_{\ovbff s}$ an $\zg$-module
$$M=\bigoplus_{\ovbff s\in\Gamma}U_{\ovbff s}\otimes\tone s\ZZ$$
with the $\g$-action
\begin{equation}\label{eq5.4}
  \begin{aligned}
    &\lone m(v_{\ovbff s}\otimes\ttwo sn)=\left((\gm\mid\al+\bff n+\bff s)\mathrm{Id}+
            \sum_{\bff p\in\ndz}\frac{\bff m^{\bff p}}{\bff p!}\rho(\xonedgm p)\right)
            v_{\ovbff s}\otimes\tthree snm;\\
    &\ltwo mr(v_{\ovbff s}\otimes\ttwo sn)=\sum_{\bff p\in\nd}
            \frac{\bff m^{\bff p}}{\bff p!}\rho(\xone p\tonebar r)v_{\ovbff s}\otimes\tfour snmr,
  \end{aligned}
\end{equation}
where $\bff m,\bff n\inr, \bff r\notinr,\bff s\in\Gamma_0$ and $v_{\ovbff s}\in U_{\ovbff s}$.
\end{thm}

Clearly, the subset of $\bff s\in\Gamma_0$ with $U_{\ovbff s}\neq0$
equals to the subset of $\bff s\in\Gamma_0$ with $M_{\al+\bff s}\neq0$.
Furthermore, we have
\begin{thm}\label{thm5.4}
Every irreducible cuspidal $\zg$-module is isomorphic to some $\vabw$,
where $\al\in\cd,\be\in\C$ and $W$ is a finite dimensional $\Gamma$-graded
irreducible $\gln$-module.
\end{thm}
\pf{
By Theorem \ref{thm5.3}, $U=\bigoplus_{\ovbff s\in\Gamma}U_{\ovbff s}$
is a finite dimensional $\Gamma$-graded irreducible $\call$-module.
Hence $\call_+U=0$ by Theorem \ref{thm4.1}(5).
Notice that $\call_+$ is also $\Gamma$-graded.
It reduces $U$ to a $\Gamma$-graded irreducible module over $\call/\call_+$,
which is isomorphic to $\gldgm\oplus\gln$.

On the other hand, since $\callt$ is an ideal of $\call$,
we obtain $\callt U=0$ or $U$ by the irreducibility of $U$.
If $\callt U=0$, then $\lone sM=0$ for all $\bff s\notinr$.
Hence $M$ is a $\gr$-module,
and the result follows from \cite{BF2}.

Suppose $\callt U=U$. 
Notice that the ideal $\call_+$ is $\Gamma$-graded and annihilates $U$.
Therefore, $U$ becomes a $\Gamma$-graded module over $\callt/\call_+$,
which is isomorphic to $\gln$. 
Then the equation (\ref{eq5.4}) implies
$\ltwo mr(v_{\ovbff s}\otimes\ttwo sn)=\bxone r v_{\ovbff s}\otimes\tfour snmr$.

We identify elements of $\gldgm$ and that of $\callx/\call_+$ through
the isomorphism given in (\ref{eq4.2}).
Notice that $\gldgm$ is solvable.
By Lie's Theorem, there exists a common eigenvector $v\in U$ for $\gldgm$ such that
$$(\e_i\gm^T)v=(x_id_{\gm})v=\be_iv,\ \be_i\in\C,\ \text{ for all }1\leq i\leq d.$$
Set $\be=\frac{\be_1}{\gm_1}$. Since
$0=[x_id_{\gm},x_1d_{\gm}]v=(\be_i\gm_1-\be_1\gm_i)v$,
we obtain $\be_i=\gm_i\be$ for all $1\leq i\leq d$.
Since $U$ is irreducible as an $\call/\call_+$-module,
$v$ can generate any vector in $U$ only by applying $\tonebar s$ and $\e_i\gm^T$.
For any $\ovbff s\in\Gamma$, we have
$$(\e_i\gm^T)(\tonebar s v)=(x_id_{\gm})(\tonebar s v)
  =[x_id_{\gm},\tonebar s]v+\tonebar s (x_id_{\gm})v
  =(\gm\mid\bff s)(x_i\tonebar s)v+\be\gm_i\tonebar sv=\be\gm_i\tonebar sv.
$$
This implies that all vectors in $U$ are common eigenvectors for $\gldgm$,
and that
$$(\bff m\gm^T)w=(\sum_{i=1}^dm_ix_id_{\gm})w=\be(\gm\mid\bff m)w,\ \
\text{ for any }w\in U,\bff m\in\zd.$$
It remains only to show that $U$ is irreducible as a $\Gamma$-graded $\gln$-module.
Let $V$ be a nonzero $\Gamma$-graded $\gln$-submodule.
Since all vectors in $V$ are common eigenvectors for $\gldgm$,
$V$ becomes a $\gldgm\oplus\gln$-module,
hence equals to $U$ by the irreducibility of $U$.
This completes the proof.
}

\section{Classification of irreducible cuspidal $\g$-modules}
\def\theequation{6.\arabic{equation}}
\setcounter{equation}{0}

In this last section we use the modified cover method to prove the following
\begin{thm}\label{thm6.1}
Let $M$ be an irreducible cuspidal $\g$-module.
Then $M$ is isomorphic to an irreducible subquotient of some $\vabw$,
where $\al\in\cd,\be\in\C$ and $W$ is a finite dimensional $\Gamma$-graded
irreducible $\gln$-module.
\end{thm}

From now on we fix an cuspidal $\g$-module $M$(not necessarily irreducible).
Recall the notion of $\zg$-module.
The proof of the following Lemma is just a elementary check, and we omit it.

\begin{lem}\label{lem6.2}
 (1) The tensor product $\gr'\otimes M$ of $\g$-modules admits a $\zg$-module structure
     if we define the $\ZZ$-action by
     $$\tone n(\lone s\otimes v)=\ltwo ns\otimes v\text{ for }\bff n\inr,\bff s\notinr, v\in M.$$
 (2) The map $\pi:\gr'\otimes M\longrightarrow M$ defined by $y\otimes v\mapsto y\cdot v$
     is a $\g$-module homomorphism, and it is surjective if $\gr'M=M$.\\
 (3) The subspace $J$ of $\gr'\otimes M$ spanned by vectors with the form of a finite sum
     $\sum\lone s\otimes v_{\bff s}$, where $\bff s\notinr,v_{\bff s}\in M$,
        and $\sum\ltwo nsv_{\bff s}=0\text{ for all }\bff n\inr$,
     is a $\zg$-submodule of $\gr'\otimes M$, and lies in $\ker\pi$.
\end{lem}

Denote $\hat M=(\gr'\otimes M)/J$.
We call this quotient $\zg$-module the $\ZZ$-cover of $M$.

\begin{prop}\label{prop6.3}
The $\zg$-module $\hat M$ is cuspidal.
\end{prop}
\pf{
Notice that $\Gamma$ is a finite group.
Consider $M$ as a cuspidal $\gr$-module.
Then $\hat M$ is a cuspidal $\ZZ\gr$-module,
and this proposition follows from \cite{BF2}.
}

\noindent
{\bf Proof of Theorem \ref{thm6.1}}:
Now assume further that $M$ is irreducible.
If $\gr'M=0$, then $M$ is irreducible as a $\gr$-module.
The result follows from \cite{BF2} or \cite{Su}.

Assume $\gr'M\neq0$. Then $\gr'M=M$ by the irreducibility of $M$.
Since $J\subseteq\ker\pi$, we get a surjective $\g$-module homomorphism
$\hat\pi: \hat M\longrightarrow M$.
Consider the composition series for the $\zg$-module $\hat M$,
$$0=\hat M_0\subset\hat M_1\subset\cdots\subset\hat M_l=\hat M,$$
where all $\hat M_i/\hat M_{i-1}$ are irreducible $\zg$-modules.
Let $k$ be the smallest integer such that $\hat\pi(\hat M_k)\neq0$.
Then we have $\hat\pi(\hat M_k)=M$ and $\hat\pi(\hat M_{k-1})=0$,
which reduces $\hat\pi$ to a $\g$-module epimorphism from $\hat M_k/\hat M_{k-1}$ to $M$.
Then Theorem \ref{thm6.1} follows from Theorem \ref{thm5.4}.

\vspace{2cm}

\noindent
{\bf Acknowledgement}:
The author is supported by the National Natural Science Foundation of
China(No. 11626157, 11801375),
the Science and Technology Foundation of Education Department of
Jiangxi Province(No. GJJ161044).

\end{document}